\documentclass[11pt]{article}
\usepackage{amssymb,amsfonts,amsthm}
\usepackage{amsmath,amsthm,amssymb}
\usepackage{amsfonts}
\newcommand{\la}{\langle}
\newcommand{\ra}{\rangle}

\newcommand{\bbb}{{\cal B}}
\newcommand{\vk}{{van Kampen}\ }

\newtheorem{theorem}{Theorem}
\newtheorem{lemma}[theorem]{Lemma}

\theoremstyle{definition}

\newtheorem{remark}[theorem]{Remark}

\newcommand{\gr}{{\cal G}}

\newcommand{\hhh}{{\cal H}}

\newcommand{\aaa}{{\cal A}}

\newcommand{\bb}{{\cal B}}

\newcommand{\iv}{^{-1}}

\newcommand{\too}{\to }

\newcommand{\Lab}{\phi}

\begin{document}

\title{Subgroups of finitely presented groups with solvable conjugacy
problem}
\author{A.Yu.Olshanskii and M.V. Sapir\thanks{Both authors were supported in part by the NSF grant DMS 0072307.
In addition, the research of the first author was supported in
part by the Russian Fund for Basic Research 02-01-00170,  the
research of the second author was supported in part by the NSF
grant DMS 9978802 and the US-Israeli BSF grant 1999298.}}
\date{} \maketitle

\begin{abstract} We prove that every countable group with solvable power
problem embeds into a finitely presented $2$-generated group with
solvable power and conjugacy problems.
\end{abstract}

\section{Introduction}

We say that a group $G$ is recursively presented if $G=\la
x_1,x_2,...\mid R\ra$ where $R$ is a recursive set of words. In
that case the set of all equalities $u=1$ where $u$ is a word in
$x_1,x_2,...,$ that are true in $G$ is recursively enumerable.
Recursive presentability follows from solvability of the word
problem. Note that in this paper, we always consider countable
groups together with their presentation: a countable group may
have a recursive presentation and a non-recursive one.

Let $G$ be a recursively presented (but not necessarily finitely
generated) group. We shall say that $G$ has solvable {\em power
problem} if there exists an algorithm which, given $u, v$ in $G$
says if $v=u^n$ for some $n\ne 0$. Notice that solvability of power problem implies solvability of the
word problem (take $u=1$). The converse implication does not hold (see \cite{McC2} or \cite{Col2}).
Notice also that if $G$ has solvable
power problem then it has solvable {\em order problem}
that is there exists an algorithm that given $u\in G$ computes the
order of $u$ in $G$. Indeed, one can first find out if there exists an $n\ne 0$ such that
$1=u^n$. Then if such an $n$ exists, find the smallest such $n$ using an algorithm that solves the word problem
in $G$.  Notice that the
``classical" formulation of the power problem does not exclude the
case $n=0$, but that change is not significant since the case
$n=0$ is a particular case of ``our" power problem for $v=1$.

In \cite{Kour} (Problem 5.21) Collins asked whether every
torsion-free group with solvable word problem can be embedded into
a finitely presented group with solvable conjugacy problem.

In this paper, we shall give positive answer to Collins' question
under the stronger assumption of solvability of the power problem.
Adding that restriction allowed us to drop the ``torsion-free"
restriction from Collins' problem.

\begin{theorem}\label{th1} Every countable group with solvable
power problem is embeddable into a $2$-generated finitely
presented group with solvable conjugacy and power problems.
\end{theorem}

Thus embedding a group into a finitely presented group can
dramatically improve its algorithmic properties.

\begin{remark}
Notice that solvability of power problem cannot
be replaced in Theorem \ref{th1} by solvability of word problem.
Indeed, there exists an example of a group with solvable word problem that
cannot be embedded into a group with solvable conjugacy problem (this example is attributed by Collins to
Macintyre in \cite[Problem 5.21]{Kour}).
\end{remark}

We also prove the following theorem that generalizes the main
result of \cite{OScol} and gives a positive answer to Collins'
problem 5.22 from \cite{Kour}.

\begin{theorem} \label{th2} Every (countable) group with solvable conjugacy
problem can be embedded into a $2$-generated finitely presented group with
solvable conjugacy problem.
\end{theorem}

In \cite{OScol}, we proved Theorem \ref{th2} only for finitely
generated groups.

\begin{remark} One can try to prove that every countable torsion-free group $G$ with solvable word problem is embeddable into
a finitely presented group with solvable word problem (and solve
Collins' problem 5.21 from \cite{Kour}) as follows. First embed
$G$ into a group $G'$ where all non-trivial elements are conjugate
using HNN extensions as in \cite{LS}. Then use Theorem \ref{th2}
to embed $G'$ into a finitely presented group with solvable
conjugacy problem. Unfortunately this idea does not work: the
group $G'$ would not necessarily have solvable word problem.
Indeed, an HNN extension of a group has solvable word problem only
if the group has solvable membership problems for the associated
subgroups. This is the reason why we cannot avoid the solvability
of power problem in Theorem \ref{th1}.
\end{remark}
\medskip

\section{Proofs}

As usual, we are going to use \vk diagrams to represent deduction
of relations in groups. Throughout the paper, for every \vk
diagram $\Delta$, $\partial\Delta$ denotes its boundary, and for
every path $p$ in a \vk diagram, $\phi(p)$ denotes its label.

The following lemma is proved in \cite{Col1}.

\begin{lemma}\label{lm1} Let $G$ be a recursively presented group
with solvable power problem. Let $a, b$ be two elements
in $G$ of the same order. Then the HNN extension $G_{a,b}=\la
G,t\mid t\iv at=b\ra$ has solvable power problem.
\end{lemma}

%

In order to embed countable groups into 2-generated groups we use
a set of positive words in the alphabet $\{a,b\}$ which is similar
to the sets used for similar purposes in \cite{McC1}, \cite{LS}
and \cite{O1}:
\begin{equation} \label{Ai}
A_i = a^{100}b^ia^{101}b^i\dots a^{199}b^i,\;\;\;\; i=1,2,\dots .
\end{equation}
We denote by $H$ the subgroup generated by these words in the free
group $F(a,b)$. The reduced words in $F(a,b)$ representing
elements of $H$ are called $H$-{\em words}. A {\em cyclic}
$H$-{\em word} is a cyclically reduced word that is freely
conjugate to an $H$-word in $F(a,b)$. Set $\lambda=\frac1{11}$ .

\begin{lemma} \label{02}

(1) Let $UV'$ and $UV''$ be two distinct cyclic $H$-words. Then
either $$|U|<\lambda \min(|UV'|, |UV''|)$$ or the word
$V'(V'')^{-1}$ is a free conjugate in $F(a,b)$ of an $H$-word.

(2) The set of (cyclic) $H$-words and the set of their subwords
 are recursive.

 (3) Suppose that a cyclic $H$-word $W$ has prefixes $w_1$ and $w_2w_1$
 for some $|w_1|\ge\lambda |W|$. Then $w_2$ is a cyclic $H$-word.

\proof (1) This statement was proved in \cite{O1} for $\lambda
=\frac17$. The same proof works for $\lambda=\frac1{11}$ since the
factor $\frac{1}{30}$ can be replaced by $\frac{1}{50}$ in Lemma 1
of \cite{O1}.

(2) This follows from the explicit forms (\ref{Ai}) for the
generators of the subgroup $H$.

(3) This also follows from the definition (\ref{Ai}) and the small
cancellation property (1).
\endproof


\end{lemma}

Consider the following construction (cf. \cite{O1}) of an
embedding of countable groups into finitely generated groups. Let
$G=\la x_1,x_2,...\mid R\ra$ be a group. Without loss of
generality we shall assume that $R$ consists of all non-empty
relations of the group $G$. Denote by $\bar R$ be the set of words
in the alphabet $\{a,b\}$ obtained by substituting $A_i$ for $x_i$
in every word from $R$. We shall denote the group $\la a,b \mid
\bar R\ra$ by $\bar G$.

We shall prove (Lemma \ref{n1} below) that $G$ embeds into $\bar G$ and that this embedding preserves solvability of
power and conjugacy problems (Lemmas \ref{n3} and
\ref{lm2} below). Note that in the literature,
there exist embeddings of countable groups into finitely generated groups which preserve solvability of
either power problem (\cite{McC1}) or conjugacy problem (\cite{Col1}). We need to preserve
solvability of both power and conjugacy problems so formally we cannot use embeddings
from either \cite{McC1} or
\cite{Col1}. Besides, our construction is easier and it yields  2-generated groups while
constructions from \cite{McC1} and \cite{Col1}
give $3$- and $4$-generated groups respectively.

The following statement is obvious.

\begin{lemma}\label{fp} The group $\bar G$ is finitely presented provided $G$ is finitely presented.
\end{lemma}


A (disc or annular) \vk diagram over $\bar R$ will be called {\em
minimal} if it contains the minimal possible number of cells among
all diagrams with the same boundary labels.

For every $\alpha>0$, a cell $\pi$ in a disc or annular diagram
$\Delta$ is called a {\em Greendlinger $\alpha$-cell} if
$\partial\pi$ contains a subpath $p$ with $|p|\ge
\alpha|\partial\pi|$, and $p$ is a subpath of a boundary component
of $\Delta$. The path $p$ will be called a {\em Greendlinger
$\alpha$-path} of $\pi$.


We say that a diagram or a map $\Delta$ satisfies the small
cancellation condition $C'(\lambda)$ if for every two cells $\pi,
\pi'$ in $\Delta$ (possibly $\pi=\pi'$), and every common subpath
$p$ of $\partial\pi$ and $(\partial\pi')\iv$, we have
$|p|<\lambda\min(|\partial\pi|, |\partial\pi'|)$.

If $\pi=\pi'$, $p$ is a common subpath of $\partial\pi$ and
$(\partial\pi)\iv$, $|p|\ge\lambda\min(|\partial\pi|,
|\partial\pi'|)$, and $\partial\pi=q'pq''p^{-1}$, where $q'$
surrounds the hole of the annular diagram $\Delta$, then $\pi$ is
said to be a {\em hoop}.

\begin{figure}[!ht]
\centering

\unitlength .7mm 
\linethickness{0.4pt}
\ifx\plotpoint\undefined\newsavebox{\plotpoint}\fi 
\begin{picture}(145.25,85)(34,0)
\put(106.5,46.13){\oval(119.5,77.75)[]}
\put(118.31,43.5){\line(0,1){.682}}
\put(118.3,44.18){\line(0,1){.68}}
\put(118.24,44.86){\line(0,1){.676}}
\multiput(118.16,45.54)(-.0304,.1678){4}{\line(0,1){.1678}}
\multiput(118.03,46.21)(-.03118,.13281){5}{\line(0,1){.13281}}
\multiput(117.88,46.87)(-.03161,.1092){6}{\line(0,1){.1092}}
\multiput(117.69,47.53)(-.03186,.09209){7}{\line(0,1){.09209}}
\multiput(117.47,48.17)(-.03196,.07904){8}{\line(0,1){.07904}}
\multiput(117.21,48.81)(-.03197,.06871){9}{\line(0,1){.06871}}
\multiput(116.92,49.42)(-.03191,.06029){10}{\line(0,1){.06029}}
\multiput(116.6,50.03)(-.03177,.05325){11}{\line(0,1){.05325}}
\multiput(116.25,50.61)(-.03159,.04726){12}{\line(0,1){.04726}}
\multiput(115.88,51.18)(-.031354,.04207){13}{\line(0,1){.04207}}
\multiput(115.47,51.73)(-.033467,.040409){13}{\line(0,1){.040409}}
\multiput(115.03,52.25)(-.032958,.035881){14}{\line(0,1){.035881}}
\multiput(114.57,52.75)(-.032435,.03187){15}{\line(-1,0){.032435}}
\multiput(114.08,53.23)(-.036455,.032322){14}{\line(-1,0){.036455}}
\multiput(113.57,53.68)(-.040991,.032751){13}{\line(-1,0){.040991}}
\multiput(113.04,54.11)(-.04617,.03316){12}{\line(-1,0){.04617}}
\multiput(112.49,54.51)(-.05215,.03355){11}{\line(-1,0){.05215}}
\multiput(111.91,54.88)(-.0538,.03083){11}{\line(-1,0){.0538}}
\multiput(111.32,55.22)(-.06084,.03084){10}{\line(-1,0){.06084}}
\multiput(110.71,55.53)(-.06926,.03076){9}{\line(-1,0){.06926}}
\multiput(110.09,55.8)(-.07959,.03057){8}{\line(-1,0){.07959}}
\multiput(109.45,56.05)(-.09263,.03023){7}{\line(-1,0){.09263}}
\multiput(108.81,56.26)(-.10974,.02969){6}{\line(-1,0){.10974}}
\multiput(108.15,56.44)(-.13333,.02884){5}{\line(-1,0){.13333}}
\put(107.48,56.58){\line(-1,0){.673}}
\put(106.81,56.69){\line(-1,0){.678}}
\put(106.13,56.77){\line(-1,0){.681}}
\put(105.45,56.81){\line(-1,0){.682}}
\put(104.77,56.81){\line(-1,0){.681}}
\put(104.08,56.78){\line(-1,0){.679}}
\put(103.41,56.72){\line(-1,0){.675}}
\multiput(102.73,56.62)(-.1672,-.0334){4}{\line(-1,0){.1672}}
\multiput(102.06,56.49)(-.13224,-.03351){5}{\line(-1,0){.13224}}
\multiput(101.4,56.32)(-.10862,-.03353){6}{\line(-1,0){.10862}}
\multiput(100.75,56.12)(-.09151,-.03347){7}{\line(-1,0){.09151}}
\multiput(100.11,55.88)(-.07847,-.03335){8}{\line(-1,0){.07847}}
\multiput(99.48,55.62)(-.06814,-.03318){9}{\line(-1,0){.06814}}
\multiput(98.87,55.32)(-.05972,-.03296){10}{\line(-1,0){.05972}}
\multiput(98.27,54.99)(-.05268,-.03271){11}{\line(-1,0){.05268}}
\multiput(97.69,54.63)(-.04669,-.03241){12}{\line(-1,0){.04669}}
\multiput(97.13,54.24)(-.041512,-.032089){13}{\line(-1,0){.041512}}
\multiput(96.59,53.82)(-.03697,-.031732){14}{\line(-1,0){.03697}}
\multiput(96.07,53.38)(-.035296,-.033584){14}{\line(-1,0){.035296}}
\multiput(95.58,52.91)(-.03353,-.035347){14}{\line(0,-1){.035347}}
\multiput(95.11,52.41)(-.031676,-.037018){14}{\line(0,-1){.037018}}
\multiput(94.67,51.89)(-.032026,-.04156){13}{\line(0,-1){.04156}}
\multiput(94.25,51.35)(-.03234,-.04674){12}{\line(0,-1){.04674}}
\multiput(93.86,50.79)(-.03263,-.05273){11}{\line(0,-1){.05273}}
\multiput(93.5,50.21)(-.03287,-.05977){10}{\line(0,-1){.05977}}
\multiput(93.17,49.61)(-.03307,-.06819){9}{\line(0,-1){.06819}}
\multiput(92.88,49)(-.03323,-.07852){8}{\line(0,-1){.07852}}
\multiput(92.61,48.37)(-.03333,-.09156){7}{\line(0,-1){.09156}}
\multiput(92.38,47.73)(-.03336,-.10867){6}{\line(0,-1){.10867}}
\multiput(92.18,47.08)(-.03331,-.13229){5}{\line(0,-1){.13229}}
\multiput(92.01,46.42)(-.0331,-.1673){4}{\line(0,-1){.1673}}
\put(91.88,45.75){\line(0,-1){.675}}
\put(91.78,45.07){\line(0,-1){.679}}
\put(91.72,44.4){\line(0,-1){1.364}}
\put(91.69,43.03){\line(0,-1){.681}}
\put(91.74,42.35){\line(0,-1){.678}}
\put(91.81,41.67){\line(0,-1){.673}}
\multiput(91.92,41)(.02904,-.13329){5}{\line(0,-1){.13329}}
\multiput(92.07,40.33)(.02986,-.10969){6}{\line(0,-1){.10969}}
\multiput(92.25,39.68)(.03037,-.09259){7}{\line(0,-1){.09259}}
\multiput(92.46,39.03)(.03069,-.07955){8}{\line(0,-1){.07955}}
\multiput(92.71,38.39)(.03087,-.06922){9}{\line(0,-1){.06922}}
\multiput(92.98,37.77)(.03093,-.06079){10}{\line(0,-1){.06079}}
\multiput(93.29,37.16)(.03091,-.05375){11}{\line(0,-1){.05375}}
\multiput(93.63,36.57)(.03363,-.0521){11}{\line(0,-1){.0521}}
\multiput(94,36)(.03323,-.04612){12}{\line(0,-1){.04612}}
\multiput(94.4,35.44)(.032814,-.040941){13}{\line(0,-1){.040941}}
\multiput(94.83,34.91)(.032377,-.036406){14}{\line(0,-1){.036406}}
\multiput(95.28,34.4)(.031919,-.032387){15}{\line(0,-1){.032387}}
\multiput(95.76,33.91)(.035931,-.032903){14}{\line(1,0){.035931}}
\multiput(96.26,33.45)(.040459,-.033406){13}{\line(1,0){.040459}}
\multiput(96.79,33.02)(.042117,-.03129){13}{\line(1,0){.042117}}
\multiput(97.34,32.61)(.0473,-.03152){12}{\line(1,0){.0473}}
\multiput(97.9,32.23)(.0533,-.03169){11}{\line(1,0){.0533}}
\multiput(98.49,31.89)(.06033,-.03181){10}{\line(1,0){.06033}}
\multiput(99.09,31.57)(.06876,-.03187){9}{\line(1,0){.06876}}
\multiput(99.71,31.28)(.07909,-.03184){8}{\line(1,0){.07909}}
\multiput(100.35,31.03)(.09213,-.03172){7}{\line(1,0){.09213}}
\multiput(100.99,30.8)(.10924,-.03145){6}{\line(1,0){.10924}}
\multiput(101.65,30.62)(.13285,-.03097){5}{\line(1,0){.13285}}
\multiput(102.31,30.46)(.1678,-.0302){4}{\line(1,0){.1678}}
\put(102.98,30.34){\line(1,0){.677}}
\put(103.66,30.25){\line(1,0){.68}}
\put(104.34,30.2){\line(1,0){1.364}}
\put(105.7,30.2){\line(1,0){.68}}
\put(106.38,30.26){\line(1,0){.676}}
\multiput(107.06,30.35)(.1677,.0307){4}{\line(1,0){.1677}}
\multiput(107.73,30.47)(.13276,.03138){5}{\line(1,0){.13276}}
\multiput(108.39,30.63)(.10915,.03178){6}{\line(1,0){.10915}}
\multiput(109.05,30.82)(.09204,.032){7}{\line(1,0){.09204}}
\multiput(109.69,31.04)(.07899,.03208){8}{\line(1,0){.07899}}
\multiput(110.32,31.3)(.06866,.03208){9}{\line(1,0){.06866}}
\multiput(110.94,31.59)(.06024,.032){10}{\line(1,0){.06024}}
\multiput(111.54,31.91)(.0532,.03186){11}{\line(1,0){.0532}}
\multiput(112.13,32.26)(.04721,.03166){12}{\line(1,0){.04721}}
\multiput(112.7,32.64)(.042022,.031418){13}{\line(1,0){.042022}}
\multiput(113.24,33.04)(.040358,.033529){13}{\line(1,0){.040358}}
\multiput(113.77,33.48)(.035831,.033012){14}{\line(1,0){.035831}}
\multiput(114.27,33.94)(.03182,.032484){15}{\line(0,1){.032484}}
\multiput(114.75,34.43)(.032266,.036505){14}{\line(0,1){.036505}}
\multiput(115.2,34.94)(.032689,.041041){13}{\line(0,1){.041041}}
\multiput(115.62,35.47)(.03309,.04622){12}{\line(0,1){.04622}}
\multiput(116.02,36.03)(.03347,.0522){11}{\line(0,1){.0522}}
\multiput(116.39,36.6)(.03075,.05385){11}{\line(0,1){.05385}}
\multiput(116.73,37.2)(.03075,.06088){10}{\line(0,1){.06088}}
\multiput(117.03,37.8)(.03065,.06931){9}{\line(0,1){.06931}}
\multiput(117.31,38.43)(.03045,.07964){8}{\line(0,1){.07964}}
\multiput(117.55,39.07)(.03009,.09268){7}{\line(0,1){.09268}}
\multiput(117.76,39.71)(.02952,.10978){6}{\line(0,1){.10978}}
\multiput(117.94,40.37)(.02863,.13338){5}{\line(0,1){.13338}}
\put(118.08,41.04){\line(0,1){.673}}
\put(118.19,41.71){\line(0,1){.678}}
\put(118.27,42.39){\line(0,1){1.109}}
\put(134.41,43.25){\line(0,1){1.193}}
\put(134.38,44.44){\line(0,1){1.191}}
\multiput(134.31,45.63)(-.0302,.2968){4}{\line(0,1){.2968}}
\multiput(134.19,46.82)(-.02816,.1969){6}{\line(0,1){.1969}}
\multiput(134.02,48)(-.03096,.16765){7}{\line(0,1){.16765}}
\multiput(133.8,49.18)(-.03302,.14548){8}{\line(0,1){.14548}}
\multiput(133.54,50.34)(-.03112,.11521){10}{\line(0,1){.11521}}
\multiput(133.23,51.49)(-.03252,.10351){11}{\line(0,1){.10351}}
\multiput(132.87,52.63)(-.03363,.09359){12}{\line(0,1){.09359}}
\multiput(132.47,53.75)(-.032059,.078987){14}{\line(0,1){.078987}}
\multiput(132.02,54.86)(-.032888,.072447){15}{\line(0,1){.072447}}
\multiput(131.53,55.95)(-.033563,.066612){16}{\line(0,1){.066612}}
\multiput(130.99,57.01)(-.032211,.057951){18}{\line(0,1){.057951}}
\multiput(130.41,58.06)(-.032718,.053618){19}{\line(0,1){.053618}}
\multiput(129.79,59.07)(-.033124,.049634){20}{\line(0,1){.049634}}
\multiput(129.12,60.07)(-.033438,.045951){21}{\line(0,1){.045951}}
\multiput(128.42,61.03)(-.033672,.042532){22}{\line(0,1){.042532}}
\multiput(127.68,61.97)(-.032422,.037703){24}{\line(0,1){.037703}}
\multiput(126.9,62.87)(-.032568,.034902){25}{\line(0,1){.034902}}
\multiput(126.09,63.75)(-.033958,.033552){25}{\line(-1,0){.033958}}
\multiput(125.24,64.58)(-.036761,.033486){24}{\line(-1,0){.036761}}
\multiput(124.36,65.39)(-.039746,.033357){23}{\line(-1,0){.039746}}
\multiput(123.44,66.15)(-.042933,.033158){22}{\line(-1,0){.042933}}
\multiput(122.5,66.88)(-.04635,.032884){21}{\line(-1,0){.04635}}
\multiput(121.53,67.57)(-.050028,.032525){20}{\line(-1,0){.050028}}
\multiput(120.53,68.23)(-.054007,.032072){19}{\line(-1,0){.054007}}
\multiput(119.5,68.83)(-.061766,.033366){17}{\line(-1,0){.061766}}
\multiput(118.45,69.4)(-.06701,.03276){16}{\line(-1,0){.06701}}
\multiput(117.38,69.93)(-.072837,.032015){15}{\line(-1,0){.072837}}
\multiput(116.28,70.41)(-.085472,.0335){13}{\line(-1,0){.085472}}
\multiput(115.17,70.84)(-.09399,.0325){12}{\line(-1,0){.09399}}
\multiput(114.05,71.23)(-.10389,.03127){11}{\line(-1,0){.10389}}
\multiput(112.9,71.58)(-.12842,.03304){9}{\line(-1,0){.12842}}
\multiput(111.75,71.87)(-.14586,.03127){8}{\line(-1,0){.14586}}
\multiput(110.58,72.12)(-.16801,.02895){7}{\line(-1,0){.16801}}
\multiput(109.4,72.33)(-.23667,.03095){5}{\line(-1,0){.23667}}
\put(108.22,72.48){\line(-1,0){1.189}}
\put(107.03,72.59){\line(-1,0){1.192}}
\put(105.84,72.65){\line(-1,0){1.193}}
\put(104.65,72.66){\line(-1,0){1.193}}
\put(103.45,72.62){\line(-1,0){1.19}}
\multiput(102.26,72.53)(-.23715,-.02703){5}{\line(-1,0){.23715}}
\multiput(101.08,72.39)(-.19655,-.03052){6}{\line(-1,0){.19655}}
\multiput(99.9,72.21)(-.16727,-.03298){7}{\line(-1,0){.16727}}
\multiput(98.73,71.98)(-.12895,-.03091){9}{\line(-1,0){.12895}}
\multiput(97.57,71.7)(-.11483,-.0325){10}{\line(-1,0){.11483}}
\multiput(96.42,71.38)(-.09452,-.03095){12}{\line(-1,0){.09452}}
\multiput(95.28,71.01)(-.086014,-.032081){13}{\line(-1,0){.086014}}
\multiput(94.17,70.59)(-.078596,-.033006){14}{\line(-1,0){.078596}}
\multiput(93.07,70.13)(-.067543,-.031646){16}{\line(-1,0){.067543}}
\multiput(91.99,69.62)(-.062309,-.03234){17}{\line(-1,0){.062309}}
\multiput(90.93,69.07)(-.05756,-.032906){18}{\line(-1,0){.05756}}
\multiput(89.89,68.48)(-.053221,-.033361){19}{\line(-1,0){.053221}}
\multiput(88.88,67.85)(-.049232,-.033718){20}{\line(-1,0){.049232}}
\multiput(87.89,67.17)(-.043476,-.032443){22}{\line(-1,0){.043476}}
\multiput(86.94,66.46)(-.040292,-.032694){23}{\line(-1,0){.040292}}
\multiput(86.01,65.7)(-.03731,-.032873){24}{\line(-1,0){.03731}}
\multiput(85.12,64.92)(-.034508,-.032985){25}{\line(-1,0){.034508}}
\multiput(84.25,64.09)(-.033141,-.034358){25}{\line(0,-1){.034358}}
\multiput(83.42,63.23)(-.033042,-.037161){24}{\line(0,-1){.037161}}
\multiput(82.63,62.34)(-.032876,-.040144){23}{\line(0,-1){.040144}}
\multiput(81.88,61.42)(-.03264,-.043329){22}{\line(0,-1){.043329}}
\multiput(81.16,60.46)(-.032324,-.046742){21}{\line(0,-1){.046742}}
\multiput(80.48,59.48)(-.033601,-.053069){19}{\line(0,-1){.053069}}
\multiput(79.84,58.47)(-.033166,-.05741){18}{\line(0,-1){.05741}}
\multiput(79.24,57.44)(-.032621,-.062162){17}{\line(0,-1){.062162}}
\multiput(78.69,56.38)(-.031952,-.067399){16}{\line(0,-1){.067399}}
\multiput(78.18,55.31)(-.033361,-.078446){14}{\line(0,-1){.078446}}
\multiput(77.71,54.21)(-.03247,-.085868){13}{\line(0,-1){.085868}}
\multiput(77.29,53.09)(-.03137,-.09437){12}{\line(0,-1){.09437}}
\multiput(76.91,51.96)(-.03302,-.11468){10}{\line(0,-1){.11468}}
\multiput(76.58,50.81)(-.03149,-.12881){9}{\line(0,-1){.12881}}
\multiput(76.3,49.65)(-.03373,-.16712){7}{\line(0,-1){.16712}}
\multiput(76.06,48.48)(-.03141,-.19641){6}{\line(0,-1){.19641}}
\multiput(75.87,47.3)(-.0281,-.23703){5}{\line(0,-1){.23703}}
\put(75.73,46.12){\line(0,-1){1.19}}
\put(75.64,44.93){\line(0,-1){2.386}}
\put(75.6,42.54){\line(0,-1){1.192}}
\put(75.65,41.35){\line(0,-1){1.189}}
\multiput(75.75,40.16)(.02988,-.23681){5}{\line(0,-1){.23681}}
\multiput(75.9,38.98)(.03288,-.19617){6}{\line(0,-1){.19617}}
\multiput(76.1,37.8)(.03061,-.146){8}{\line(0,-1){.146}}
\multiput(76.35,36.63)(.03245,-.12857){9}{\line(0,-1){.12857}}
\multiput(76.64,35.48)(.0308,-.10403){11}{\line(0,-1){.10403}}
\multiput(76.98,34.33)(.03208,-.09414){12}{\line(0,-1){.09414}}
\multiput(77.36,33.2)(.033112,-.085623){13}{\line(0,-1){.085623}}
\multiput(77.79,32.09)(.031685,-.072981){15}{\line(0,-1){.072981}}
\multiput(78.27,30.99)(.032456,-.067158){16}{\line(0,-1){.067158}}
\multiput(78.79,29.92)(.033086,-.061916){17}{\line(0,-1){.061916}}
\multiput(79.35,28.87)(.033595,-.05716){18}{\line(0,-1){.05716}}
\multiput(79.95,27.84)(.032298,-.050175){20}{\line(0,-1){.050175}}
\multiput(80.6,26.83)(.032673,-.046498){21}{\line(0,-1){.046498}}
\multiput(81.29,25.86)(.032963,-.043083){22}{\line(0,-1){.043083}}
\multiput(82.01,24.91)(.033176,-.039896){23}{\line(0,-1){.039896}}
\multiput(82.77,23.99)(.033319,-.036913){24}{\line(0,-1){.036913}}
\multiput(83.57,23.11)(.033398,-.034109){25}{\line(0,-1){.034109}}
\multiput(84.41,22.25)(.034754,-.032726){25}{\line(1,0){.034754}}
\multiput(85.28,21.44)(.037556,-.032593){24}{\line(1,0){.037556}}
\multiput(86.18,20.65)(.040536,-.032392){23}{\line(1,0){.040536}}
\multiput(87.11,19.91)(.0458,-.033646){21}{\line(1,0){.0458}}
\multiput(88.07,19.2)(.049483,-.033348){20}{\line(1,0){.049483}}
\multiput(89.06,18.53)(.053469,-.032961){19}{\line(1,0){.053469}}
\multiput(90.08,17.91)(.057805,-.032473){18}{\line(1,0){.057805}}
\multiput(91.12,17.32)(.06255,-.031872){17}{\line(1,0){.06255}}
\multiput(92.18,16.78)(.072297,-.033215){15}{\line(1,0){.072297}}
\multiput(93.27,16.28)(.078841,-.032416){14}{\line(1,0){.078841}}
\multiput(94.37,15.83)(.086252,-.031436){13}{\line(1,0){.086252}}
\multiput(95.49,15.42)(.10336,-.03298){11}{\line(1,0){.10336}}
\multiput(96.63,15.06)(.11507,-.03164){10}{\line(1,0){.11507}}
\multiput(97.78,14.74)(.14533,-.03368){8}{\line(1,0){.14533}}
\multiput(98.94,14.47)(.16751,-.03172){7}{\line(1,0){.16751}}
\multiput(100.12,14.25)(.19677,-.02905){6}{\line(1,0){.19677}}
\multiput(101.3,14.08)(.2967,-.0316){4}{\line(1,0){.2967}}
\put(102.48,13.95){\line(1,0){1.191}}
\put(103.67,13.87){\line(1,0){1.193}}
\put(104.87,13.84){\line(1,0){1.193}}
\put(106.06,13.86){\line(1,0){1.192}}
\put(107.25,13.93){\line(1,0){1.188}}
\multiput(108.44,14.04)(.23643,.03272){5}{\line(1,0){.23643}}
\multiput(109.62,14.21)(.16779,.0302){7}{\line(1,0){.16779}}
\multiput(110.8,14.42)(.14563,.03236){8}{\line(1,0){.14563}}
\multiput(111.96,14.68)(.11535,.0306){10}{\line(1,0){.11535}}
\multiput(113.11,14.98)(.10365,.03205){11}{\line(1,0){.10365}}
\multiput(114.25,15.34)(.09374,.03321){12}{\line(1,0){.09374}}
\multiput(115.38,15.74)(.079131,.031701){14}{\line(1,0){.079131}}
\multiput(116.49,16.18)(.072595,.03256){15}{\line(1,0){.072595}}
\multiput(117.58,16.67)(.066763,.033261){16}{\line(1,0){.066763}}
\multiput(118.64,17.2)(.058096,.031949){18}{\line(1,0){.058096}}
\multiput(119.69,17.77)(.053765,.032475){19}{\line(1,0){.053765}}
\multiput(120.71,18.39)(.049783,.032899){20}{\line(1,0){.049783}}
\multiput(121.71,19.05)(.046102,.03323){21}{\line(1,0){.046102}}
\multiput(122.68,19.75)(.042684,.033479){22}{\line(1,0){.042684}}
\multiput(123.61,20.48)(.039495,.033653){23}{\line(1,0){.039495}}
\multiput(124.52,21.26)(.035049,.03241){25}{\line(1,0){.035049}}
\multiput(125.4,22.07)(.033705,.033805){25}{\line(0,1){.033805}}
\multiput(126.24,22.91)(.033652,.036609){24}{\line(0,1){.036609}}
\multiput(127.05,23.79)(.033536,.039594){23}{\line(0,1){.039594}}
\multiput(127.82,24.7)(.033352,.042783){22}{\line(0,1){.042783}}
\multiput(128.55,25.64)(.033093,.046201){21}{\line(0,1){.046201}}
\multiput(129.25,26.61)(.032751,.049881){20}{\line(0,1){.049881}}
\multiput(129.9,27.61)(.032316,.053861){19}{\line(0,1){.053861}}
\multiput(130.52,28.64)(.033645,.061614){17}{\line(0,1){.061614}}
\multiput(131.09,29.68)(.033063,.066861){16}{\line(0,1){.066861}}
\multiput(131.62,30.75)(.032344,.072691){15}{\line(0,1){.072691}}
\multiput(132.11,31.84)(.031466,.079225){14}{\line(0,1){.079225}}
\multiput(132.55,32.95)(.03293,.09384){12}{\line(0,1){.09384}}
\multiput(132.94,34.08)(.03174,.10375){11}{\line(0,1){.10375}}
\multiput(133.29,35.22)(.03362,.12827){9}{\line(0,1){.12827}}
\multiput(133.59,36.37)(.03193,.14572){8}{\line(0,1){.14572}}
\multiput(133.85,37.54)(.02971,.16788){7}{\line(0,1){.16788}}
\multiput(134.06,38.71)(.03202,.23653){5}{\line(0,1){.23653}}
\put(134.22,39.9){\line(0,1){1.188}}
\put(134.33,41.09){\line(0,1){2.165}}
\put(99.5,41.25){\rule{12\unitlength}{4\unitlength}}
\put(83.5,44.75){\vector(-1,0){.07}}\put(91.5,44.75){\line(-1,0){16}}
\put(85,40){\makebox(0,0)[cc]{$p$}}
\put(83.75,63.63){\vector(-2,-3){.07}}\multiput(84.5,64.75)(-.0333333,-.05){45}{\line(0,-1){.05}}
\put(103.38,57){\vector(4,1){.07}}\multiput(102,56.75)(.183333,.033333){15}{\line(1,0){.183333}}
\put(103,53){\makebox(0,0)[cc]{$q'$}}
\put(88.25,62.5){\makebox(0,0)[cc]{$q''$}}
\put(147.5,68.25){\makebox(0,0)[cc]{$\Delta$}}
\put(117,31){\makebox(0,0)[cc]{$\pi$}}
\end{picture}

\caption{A hoop.} \label{ahoop}
\end{figure}

\begin{lemma} \label{n1} (1) A minimal disc or annular diagram
over $\bar R$, having no hoops in annular case, satisfies the
$C'(\lambda)$-condition.

(2) The labels of the subpaths $q'$ and $q''$ in the boundary of a
hoop (see the definition and notation before the lemma) are cyclic
$H$-words.

(3) If a boundary component $q$ of a diagram (annular diagram)
over ${\bar R}$ is an $H$-word (a cyclic $H$-word), and a cell
$\pi$ has boundary $p_1p_2$, where $q=p_1q'$ and
$|p_1|\ge\lambda|\partial\pi|$, then the label of the path
$p_2^{-1}q'$ is freely equal (freely conjugate) to an $H$-word.

(4) If a disc map $\Delta$ satisfies $C'(\lambda)$-condition and
contains a cell, then (a) it has a Greendlinger
$(1-3\lambda)$-cell, and (b) the number of (non-directed) edges of
$\Delta$ does not exceed \\
$(1+3\lambda(1-6\lambda)^{-1})|\partial\Delta|$.

(5) If an annular map $\Delta$  contains at least one cell and
satisfies $C'(\lambda)$-condition, then either $\Delta$ has a
Greendlinger $(1-4\lambda)$-cell or every cell $\pi$ of $\Delta$
has boundary subpaths $p_1$ and $p_2$ on both boundary components;
in particular $|p_1|+|p_2|>(1-2\lambda)|\partial\pi|$.

(6)(a) The mapping $x_i\too A_i$ extends to an embedding of the
group $G$ into $\bar G$. (b) Under this embedding, two elements of
$G$ that are conjugate in $\bar G$ are also conjugate in $G$.

\end{lemma}

\proof Lemma \ref{02} implies assertions (1), (2), (3) and (6)(a)
as this was shown in $\S 2$ of \cite{O1}. The proof of (6)(b) is
similar to the proof of (6)(a) but one should take an annular
diagram instead of a disc one and use assertion (2) when
considering annular subdiagrams between two hoops. (The statement
(6)(b) was also proved by Ilya Belyaev in \cite{Bel}.)
 The proof of assertion (4)(a) can be found in \cite{LS} (see
Theorem 4.4, formulated in terms of relations). The assertion
(4)(b) follows from (4)(a) by straightforward induction on the
number of cells in $\Delta$. The assertion (5) is contained in
\cite[Theorems 5.3 and 5.5]{LS}. (In the proofs of these theorems,
one can replace $C'(1/6)$ by $C'(\lambda)$, where $\lambda\le
1/6$. Accordingly one can replace $(1/2)R$ by $(1-3\lambda)R$ in
the assumptions of these theorems. Since $i(D)=p/q +2 =4$ in the
formulation of \cite[Theorem 5.3]{LS}, $D$ is a Greendlinger
$(1-4\lambda)$-cell.)\endproof

We call a (cyclically) reduced word $w$ {\em (cyclically) $\bar
R$-reduced} if it has no (cyclic) subword $v$, where $v$ is a
subword of a relator $r\in \bar R$ and $|v|>\frac12 |r|.$ If the
word problem is decidable for $G$, then, by Lemma \ref{02}(2), for
every word $w$, one can effectively find a (cyclically) $\bar
R$-reduced word $w'$ which is equal (or is conjugate) to $w$ in
the group $\bar G$.

\begin{lemma}\label{n2} Assume that a $\bar R$-reduced (cyclically $\bar
R$-reduced) word $w$ is equal (is conjugate) in $\bar G$ to an
$H$-word (a cyclic $H$-word). Then $w$ is equal (is conjugate) to
an $H$-word (a cyclic $H$-word) in the free group $F(a,b)$.
\end{lemma}

\proof We consider only the cyclic case. Let $\Delta$ be a minimal
diagram for the conjugation of $w$ and some cyclic $H$-word $w'$,
and let $\Delta$ have minimal number of cells over all such $w'$.
By Lemma \ref{n1}(2), $\Delta$ has no hoops. If $\Delta$ has a
cell $\pi$, then, by Lemma \ref{n1} (3) and by the minimality of
$\Delta$, $\pi$ cannot have a common boundary subpath of length at
least $\lambda|\partial\pi|$ with the contour of $\Delta$ labeled
by $w'$. Then, by Lemma \ref{n1} (2),(5), $\pi$ must have a
boundary subpath of length greater than
$(1-4\lambda)|\partial\pi|$, lying on the contour of $\Delta$
labeled by $w'$. This contradicts the cyclic $\bar R$-reducibility
of $w$, since $\lambda \le 1/8$.
\endproof

\begin{lemma}\label{n3} If the group $G$ has decidable word or conjugacy
problem, then so has the group $\bar G$.
\end{lemma}

\proof Here we consider only the conjugacy problem. Let $u$ and
$v$ be two words under our investigation. We may assume that they
are cyclically $\bar R$-reduced since the conjugacy (and the word)
problem is decidable for $G$. By Lemma \ref{n2}, we can also
assume that they are $H$-words if they are conjugate to $H$-words
in $\bar G$. In the later case, it suffices to check conjugacy of
$u$ and $v$ in $G$ by Lemma \ref{n1} (6)(b). If $u$ and $v$ are
conjugate in $\bar G$, but none of them is conjugate to a cyclic
$H$-word, then a minimal diagram $\Delta$ for the conjugation of
$u$ and $v$ has no hoops by Lemma \ref{n1} (2). It also has no
Greendlinger $\alpha$-cells for $\alpha>1/2$, since $u$ and $v$
are cyclically $\bar R$-reduced. It follows from Lemma \ref{n1},
part  (5), that the sum of perimeters of the cells of $\Delta$
does not exceed $(1-2\lambda^{-1})(|u|+|v|)$, and therefore such
diagrams can be checked by exhaustion. \endproof

\begin{lemma}\label{lm2} Let $G=\la X\mid R\ra$ be a recursively presented
group with solvable power problem. Then the group $\bar
G$ has solvable power problem.
\end{lemma}

\proof Let $u$ and $v$ be words in the alphabet $\{a,b\}$. Suppose
that $v=u^n$ in $\bar G$, where $n\ge 1$. By Lemmas  \ref{02} (2)
and \ref{n2}, we may assume that the word $u$ is cyclically $\bar
R$-reduced and it is a cyclic $H$-word if it is a conjugate of an
$H$-word in $\bar G$. If $u$ is an $H$-word, then so is $v$, and,
by Lemma \ref{n1} (6)(a), we may refer to the solvability of the
power problem in $G$. Therefore we further assume that $u$ is not
a conjugate of an $H$-word in $\bar G$.

Let $\Delta$ be a minimal diagram over $\bar R$ whose contour is
$q_1q_2$, where $q_1$ and $q_2$ are labeled by $u^n$ and $v\iv$,
respectively. Call a cell $\pi$ of $\Delta$ {\em suitable}, if its
boundary has a common subpath $p$ with $q_1$, and
$|p|>(1/2+\lambda)|\partial\pi|$.

Suppose $\Delta$ has a suitable cell $\pi$. If $|p|\ge |u|+\lambda
|\partial\pi|$, then $u$ is a cyclic $H$-word by Lemma \ref{02}
(3), since the label of $p$ is a subword of $u^n$; a
contradiction. But it follows from the inequality $|p|<|u|+\lambda
|\partial\pi|$ that $|p|<|u|+\lambda(1/2+\lambda)^{-1}|p|$, i.e.
$|p|< (1-\lambda(1/2+\lambda)^{-1})^{-1}|u|$. It also follows that
the word $u$ is not cyclically $\bar R$-reduced, since an
application to a cyclic permutation of $u$ of the $\bar R$-relator
corresponding to $\pi$, gives a word of length at most
$$|\partial\pi|-|p|+\lambda |\partial\pi| <
\frac12|\partial\pi|\le \frac12(1/2+\lambda)^{-1}|p|<
(1+2\lambda)^{-1}(1-\lambda(1/2+\lambda)^{-1})^{-1}|u|=|u|$$ The
contradiction shows that $\Delta$ has no suitable cells.

Now assume that $\Delta$ has a cell $\Pi$ having two maximal
boundary subpaths $p_1$ and $p_2$ on $q_1$ (a ``bad" cell). Then
there must be cells in the subdiagram $\Gamma$ between $\Pi$ and
$q_1$, and one may chose $\Pi$ so that there are no bad cells in
$\Gamma$. Similarly we assume that there are no cells in $\Gamma$
having two maximal boundary subpaths on $\partial\Pi$, since
otherwise we can decrease the number of cells in $\Gamma$. Then by
Lemma \ref{n1} (1), (4)(a), the diagram $\Delta$ has a suitable
cell inside $\Gamma$ since $1-3\gamma-\lambda\ge 1/2+\lambda$.

\begin{figure}[!ht]
\centering
\unitlength .6mm 
\linethickness{0.4pt}
\ifx\plotpoint\undefined\newsavebox{\plotpoint}\fi 
\begin{picture}(125.88,85.63)(14,0)
\qbezier(41.5,19.5)(29.5,42.63)(43.5,64.25)
\qbezier(43.5,64.25)(79.63,85.63)(114.25,63.5)
\qbezier(114.25,63.5)(125.88,40.38)(111,18.75)
\qbezier(111,18.75)(81.38,2.38)(41.25,19.5)
\qbezier(43.5,64)(77.63,53.5)(114.25,63)
\qbezier(41.5,19.5)(77.88,32.38)(110.75,18.75)
\put(76.25,74.75){\vector(1,0){.07}}\put(74.75,74.75){\line(1,0){3}}
\put(35.75,40.38){\vector(0,-1){.07}}\put(35.75,42.75){\line(0,-1){4.75}}
\put(80.63,11){\vector(-1,0){.07}}\put(83,11){\line(-1,0){4.75}}
\put(119.25,42){\vector(0,1){.07}}\put(119.25,40.25){\line(0,1){3.5}}
\put(29.25,41.5){\makebox(0,0)[cc]{$p_2$}}
\put(71.5,80){\makebox(0,0)[cc]{$q_1$}}
\put(123.75,41.25){\makebox(0,0)[cc]{$p_1$}}
\put(78.5,6){\makebox(0,0)[cc]{$q_2$}}
\put(77.25,43.25){\makebox(0,0)[cc]{$\Pi$}}
\put(77.75,67.5){\makebox(0,0)[cc]{$\Gamma$}}
\end{picture}

\caption{A ``bad" cell.} \label{fig6}
\end{figure}
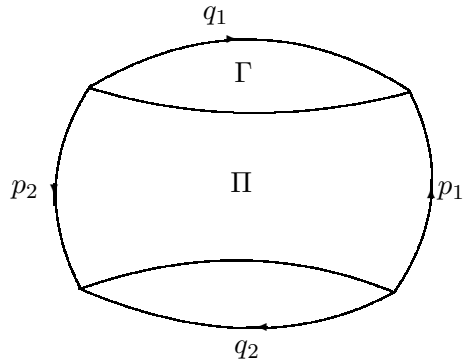

This contradiction shows that there are no bad cells in $\Delta$.
The same consideration shows that the path $q_1$ is simple. Then
$\Delta$ has at least $(1/2+\lambda)^{-1}(|q_1|-|q_2|)$ edges
because it has no suitable cells. This inequality and Lemma
\ref{n1} (4)(b) give a linear upper bound for $|q_1|$ in terms of
$|q_2|$ since $(1/2+\lambda)^{-1}>1+3\lambda(1-6\lambda)^{-1}$.
This reduces the problem to the word problem in $\bar G$ which is
decidable by Lemma \ref{n3}.
\endproof

\begin{lemma} \label{lm3} Let $G$ be a finitely generated
recursively presented group with solvable power and conjugacy
problems. Then the group $G$ can be embedded into a finitely
presented group $H$ with solvable power and conjugacy problems.
\end{lemma}

\proof Let us use the embedding from \cite{OScol}.  We are going
to use the notation and results from \cite{OScol}. In particular,
from now on we shall denote $G$ by $\gr$, and $H$ by $\hhh$ as in
\cite{OScol}. We have proved in \cite{OScol} that $\hhh$ has
solvable conjugacy problem. It remains to prove the solvability of
the power problem.

Recall that the set of generators of $\hhh$ consists of
$k$-letters, $a$-letters, $\theta$-letters and $x$-letters. The
subgroup $\gr$ in $\hhh$ is generated by a subset $\aaa(P_1)$ of
the set of $a$-letters. As in \cite{OScol}, we include all
relations of $\gr$ into the presentation of $\hhh$.

Among the relations of $\hhh$, there is one, called the {\em hub}
which is a word in $k$-letters of length $N$ (in \cite{OScol} $N$
is any even number $\ge 8$; here we take $N\ge 14$), all letters
occurring in the hub are different. Every non-hub cell in a \vk
diagram over $\hhh$ that contains an edge labeled by a $k$-letter
$K$ (i.e. a {\em $k$-edge}) also contains an edge labeled by
$K\iv$, so we can consider $k$-bands in $\Delta$ (see the precise
definition of a band in \cite{OScol}). Similarly, we can consider
$\theta$-bands and $a$-bands. The group given by the presentation
of $\hhh$ without the hub is denoted by $\hhh_1$.

1. Let $u$ be a word in generators of $\hhh$. Consider a word $u'$
that is a conjugate of $u$ in $\hhh$ and has minimal number of
$k$-letters among all words in the conjugacy class of $\hhh$.
Consider a minimal annular diagram $\Delta$ for this conjugation
with contours $p$ and $p'$, $\phi(p)=u, \phi(p')=u'$. By the
minimality in the choice of $u'$, every hub $\Pi$ of $\Delta$ has
at most $N/2$ $k$-bands starting on $\Pi$ and ending on $p'$.
Assume there exists a hub in $\Delta$. Since $N\ge 10$,
\cite[Lemmas 10.4 and 10.3]{OScol} provide us with a hub $\Pi$
connected with $p$ by two $k$-bands, such that there are no hubs
between these $k$-bands. By \cite[Lemma 10.5]{OScol}, the hub
$\Pi$ can be effectively cut out of $\Delta$ with a recursive
replacement of $u$ by a conjugate word. However, by \cite[Lemmas
10.4 and 10.3]{OScol}, the number of hubs in $\Delta$ is not
greater than the doubled number of $k$-letters in $u$. Therefore
we can recursively obtain a word $u''$, such that $u''$ and $u'$
are conjugate in the group $\hhh_1$. Now, by \cite[Lemma
5.6]{OScol}, starting with $u''$, we can recursively obtain an
$\hhh_1$-conjugate word $u'''$ having minimal number of
$k$-letters in its $\hhh$-conjugacy class and minimal number of
$\theta$-letters in its $\hhh_1$-conjugacy class. By \cite[Lemma
5.6]{OScol}, $u'''$ is also not conjugate to a word with fewer
$a$-letters if a deduction of the latest conjugation does not
employ $\theta$-relations. We shall call a word $u'''$ with these
three properties {\em cyclically minimal}. Similarly for every
word $w$ in the generators of $\hhh$, we can effectively find a
word $w'$ which is equal to $w$ in $\hhh$ and has minimal number
of $k$-letters among all words that are equal to $w$ in $\hhh$,
minimal number of $\theta$-letters among all words that are equal
to $w'$ in $\hhh_1$, and not equal in $\hhh_1$ to a word with
fewer $a$-letters if a deduction of this equality does not use
$\theta$-relations. Such words $w'$ will be called {\em minimal}.

2. Let $u, w$ be words in the generators of $\hhh$. Suppose that
$w=u^n$ in $\hhh$ for some $n\ge 0$. To prove the lemma, we need
to recursively bound $n$ in terms of $w$ and $u$. We can assume
that the word $u$ is cyclically reduced and cyclically minimal by
part 1. We can also assume that $w$ is minimal. Consider the
corresponding minimal \vk diagram $\Delta$ such that
$\partial\Delta=pq\iv$, $\Lab(p)=w$, $\Lab(q)=u^n$.

3. We set $N'=N/2+1$. Suppose that there exists a hub $\pi$ in
$\Delta$ such that some consecutive $k$-bands
$\bbb_1,...,\bbb_{N'}$ starting on $\partial\pi$ end on $q$, and
between two consecutive $k$-bands starting on $\partial\pi$, there
are no other $k$-bands. We shall call these hubs {\em $u$-close}.
In particular, it implies that between any two consecutive bands
$\bb_i$, there are no hubs (since $k$-bands do not intersect).
Since all $k$-edges on $\partial\pi$ have different labels, and
there are no other $k$-edges between the end edges of
$\bbb_1,...,\bbb_{N'}$ on $q$, we can conclude that all the bands
$\bbb_1,...,\bbb_{N'}$ connect $\partial\pi$ with a subpath of $q$
labeled by a cyclic shift of $u$, contrary to the assumption that
$u$ is minimal, because $N'>N-N'$. Hence $\Delta$ contains no
$u$-close hubs. Since $N\ge 14$, this result, \cite[Lemma
10.4]{OScol} and \cite[Lemma 3.4]{Ol97} give a linear upper bound
for the number of hubs of $\Delta$ in terms of $|w|$.

\begin{figure}[!ht]
\centering
\unitlength .7mm 
\linethickness{0.4pt}
\ifx\plotpoint\undefined\newsavebox{\plotpoint}\fi 
\begin{picture}(142.75,81.5)(0,0)
\put(73.25,14){\vector(1,0){.07}}\put(37.5,14){\line(1,0){71.5}}
\put(29.5,61.25){\vector(2,3){.07}}\qbezier(37.25,13.75)(2,77.13)(76.75,77)
\put(117.75,61.25){\vector(1,-2){.07}}\qbezier(76.75,77)(142.75,76.88)(108.75,14.25)
\put(83.36,36.25){\line(0,1){.568}}
\put(83.34,36.82){\line(0,1){.567}}
\put(83.3,37.38){\line(0,1){.563}}
\put(83.22,37.95){\line(0,1){.559}}
\multiput(83.11,38.51)(-.02716,.11041){5}{\line(0,1){.11041}}
\multiput(82.98,39.06)(-.03304,.10879){5}{\line(0,1){.10879}}
\multiput(82.81,39.6)(-.03235,.08906){6}{\line(0,1){.08906}}
\multiput(82.62,40.14)(-.03178,.07474){7}{\line(0,1){.07474}}
\multiput(82.4,40.66)(-.03127,.06381){8}{\line(0,1){.06381}}
\multiput(82.15,41.17)(-.03079,.05515){9}{\line(0,1){.05515}}
\multiput(81.87,41.67)(-.03371,.05342){9}{\line(0,1){.05342}}
\multiput(81.57,42.15)(-.03287,.04639){10}{\line(0,1){.04639}}
\multiput(81.24,42.61)(-.0321,.04051){11}{\line(0,1){.04051}}
\multiput(80.88,43.06)(-.03137,.0355){12}{\line(0,1){.0355}}
\multiput(80.51,43.48)(-.03323,.03377){12}{\line(0,1){.03377}}
\multiput(80.11,43.89)(-.03499,.03194){12}{\line(-1,0){.03499}}
\multiput(79.69,44.27)(-.03998,.03275){11}{\line(-1,0){.03998}}
\multiput(79.25,44.63)(-.04585,.03362){10}{\line(-1,0){.04585}}
\multiput(78.79,44.97)(-.04758,.03111){10}{\line(-1,0){.04758}}
\multiput(78.31,45.28)(-.05464,.03169){9}{\line(-1,0){.05464}}
\multiput(77.82,45.56)(-.0633,.0323){8}{\line(-1,0){.0633}}
\multiput(77.32,45.82)(-.07421,.03299){7}{\line(-1,0){.07421}}
\multiput(76.8,46.05)(-.07587,.02897){7}{\line(-1,0){.07587}}
\multiput(76.27,46.26)(-.0902,.029){6}{\line(-1,0){.0902}}
\multiput(75.72,46.43)(-.10995,.02895){5}{\line(-1,0){.10995}}
\put(75.18,46.58){\line(-1,0){.557}}
\put(74.62,46.69){\line(-1,0){.562}}
\put(74.06,46.78){\line(-1,0){.566}}
\put(73.49,46.83){\line(-1,0){.568}}
\put(72.92,46.86){\line(-1,0){.568}}
\put(72.35,46.85){\line(-1,0){.567}}
\put(71.79,46.81){\line(-1,0){.565}}
\put(71.22,46.75){\line(-1,0){.56}}
\multiput(70.66,46.65)(-.1385,-.0317){4}{\line(-1,0){.1385}}
\multiput(70.11,46.52)(-.10932,-.03126){5}{\line(-1,0){.10932}}
\multiput(69.56,46.37)(-.08957,-.0309){6}{\line(-1,0){.08957}}
\multiput(69.02,46.18)(-.07525,-.03056){7}{\line(-1,0){.07525}}
\multiput(68.5,45.97)(-.06431,-.03023){8}{\line(-1,0){.06431}}
\multiput(67.98,45.72)(-.0626,-.03363){8}{\line(-1,0){.0626}}
\multiput(67.48,45.46)(-.05396,-.03283){9}{\line(-1,0){.05396}}
\multiput(67,45.16)(-.04691,-.03211){10}{\line(-1,0){.04691}}
\multiput(66.53,44.84)(-.04102,-.03143){11}{\line(-1,0){.04102}}
\multiput(66.08,44.49)(-.03928,-.03359){11}{\line(-1,0){.03928}}
\multiput(65.64,44.12)(-.03431,-.03267){12}{\line(-1,0){.03431}}
\multiput(65.23,43.73)(-.03251,-.03446){12}{\line(0,-1){.03446}}
\multiput(64.84,43.32)(-.0334,-.03944){11}{\line(0,-1){.03944}}
\multiput(64.47,42.88)(-.03123,-.04117){11}{\line(0,-1){.04117}}
\multiput(64.13,42.43)(-.03188,-.04707){10}{\line(0,-1){.04707}}
\multiput(63.81,41.96)(-.03257,-.05412){9}{\line(0,-1){.05412}}
\multiput(63.52,41.47)(-.03333,-.06276){8}{\line(0,-1){.06276}}
\multiput(63.25,40.97)(-.02992,-.06446){8}{\line(0,-1){.06446}}
\multiput(63.01,40.46)(-.0302,-.07539){7}{\line(0,-1){.07539}}
\multiput(62.8,39.93)(-.03046,-.08972){6}{\line(0,-1){.08972}}
\multiput(62.62,39.39)(-.03074,-.10947){5}{\line(0,-1){.10947}}
\multiput(62.47,38.84)(-.031,-.1387){4}{\line(0,-1){.1387}}
\put(62.34,38.29){\line(0,-1){.561}}
\put(62.25,37.73){\line(0,-1){.565}}
\put(62.18,37.16){\line(0,-1){2.269}}
\put(62.23,34.89){\line(0,-1){.562}}
\put(62.32,34.33){\line(0,-1){.556}}
\multiput(62.44,33.78)(.02948,-.10981){5}{\line(0,-1){.10981}}
\multiput(62.58,33.23)(.02944,-.09006){6}{\line(0,-1){.09006}}
\multiput(62.76,32.69)(.02933,-.07573){7}{\line(0,-1){.07573}}
\multiput(62.97,32.16)(.03335,-.07405){7}{\line(0,-1){.07405}}
\multiput(63.2,31.64)(.03261,-.06314){8}{\line(0,-1){.06314}}
\multiput(63.46,31.13)(.03195,-.05449){9}{\line(0,-1){.05449}}
\multiput(63.75,30.64)(.03134,-.04743){10}{\line(0,-1){.04743}}
\multiput(64.06,30.17)(.03076,-.04153){11}{\line(0,-1){.04153}}
\multiput(64.4,29.71)(.03294,-.03982){11}{\line(0,-1){.03982}}
\multiput(64.76,29.27)(.03211,-.03483){12}{\line(0,-1){.03483}}
\multiput(65.15,28.85)(.03393,-.03306){12}{\line(1,0){.03393}}
\multiput(65.55,28.46)(.03565,-.0312){12}{\line(1,0){.03565}}
\multiput(65.98,28.08)(.04066,-.0319){11}{\line(1,0){.04066}}
\multiput(66.43,27.73)(.04654,-.03264){10}{\line(1,0){.04654}}
\multiput(66.89,27.41)(.05358,-.03345){9}{\line(1,0){.05358}}
\multiput(67.38,27.11)(.0553,-.03053){9}{\line(1,0){.0553}}
\multiput(67.87,26.83)(.06396,-.03096){8}{\line(1,0){.06396}}
\multiput(68.39,26.58)(.07489,-.03142){7}{\line(1,0){.07489}}
\multiput(68.91,26.36)(.08921,-.03192){6}{\line(1,0){.08921}}
\multiput(69.45,26.17)(.10895,-.03251){5}{\line(1,0){.10895}}
\multiput(69.99,26.01)(.1382,-.0333){4}{\line(1,0){.1382}}
\put(70.54,25.88){\line(1,0){.559}}
\put(71.1,25.77){\line(1,0){.564}}
\put(71.67,25.7){\line(1,0){.567}}
\put(72.23,25.66){\line(1,0){1.137}}
\put(73.37,25.66){\line(1,0){.566}}
\put(73.94,25.71){\line(1,0){.563}}
\put(74.5,25.79){\line(1,0){.558}}
\multiput(75.06,25.9)(.11027,.02769){5}{\line(1,0){.11027}}
\multiput(75.61,26.04)(.10863,.03356){5}{\line(1,0){.10863}}
\multiput(76.15,26.2)(.0889,.03278){6}{\line(1,0){.0889}}
\multiput(76.68,26.4)(.07458,.03214){7}{\line(1,0){.07458}}
\multiput(77.21,26.63)(.06366,.03158){8}{\line(1,0){.06366}}
\multiput(77.72,26.88)(.055,.03106){9}{\line(1,0){.055}}
\multiput(78.21,27.16)(.04793,.03057){10}{\line(1,0){.04793}}
\multiput(78.69,27.46)(.04623,.03309){10}{\line(1,0){.04623}}
\multiput(79.15,27.79)(.04035,.03229){11}{\line(1,0){.04035}}
\multiput(79.6,28.15)(.03535,.03154){12}{\line(1,0){.03535}}
\multiput(80.02,28.53)(.03361,.03339){12}{\line(1,0){.03361}}
\multiput(80.42,28.93)(.03177,.03514){12}{\line(0,1){.03514}}
\multiput(80.81,29.35)(.03256,.04014){11}{\line(0,1){.04014}}
\multiput(81.16,29.79)(.0334,.04601){10}{\line(0,1){.04601}}
\multiput(81.5,30.25)(.03088,.04773){10}{\line(0,1){.04773}}
\multiput(81.81,30.73)(.03142,.0548){9}{\line(0,1){.0548}}
\multiput(82.09,31.22)(.032,.06345){8}{\line(0,1){.06345}}
\multiput(82.35,31.73)(.03263,.07437){7}{\line(0,1){.07437}}
\multiput(82.57,32.25)(.03337,.08868){6}{\line(0,1){.08868}}
\multiput(82.77,32.78)(.02857,.09034){6}{\line(0,1){.09034}}
\multiput(82.95,33.32)(.02842,.11009){5}{\line(0,1){.11009}}
\put(83.09,33.87){\line(0,1){.557}}
\put(83.2,34.43){\line(0,1){.562}}
\put(83.28,34.99){\line(0,1){.566}}
\put(83.33,35.56){\line(0,1){.689}}
\multiput(83.25,35.5)(.092269327,.033665835){401}{\line(1,0){.092269327}}
\multiput(82.75,39.75)(.084474886,.033675799){438}{\line(1,0){.084474886}}
\multiput(78,45.5)(.033691406,.057617188){512}{\line(0,1){.057617188}}
\multiput(75.25,46.75)(.03372093,.068604651){430}{\line(0,1){.068604651}}
\multiput(62,36.5)(-.116185897,.033653846){312}{\line(-1,0){.116185897}}
\multiput(62.5,40.25)(-.099587912,.033653846){364}{\line(-1,0){.099587912}}
\put(58.25,58.25){\makebox(0,0)[cc]{$\cdots$}}
\put(118.63,38.88){\vector(-1,-3){.07}}\multiput(119,40)(-.032609,-.097826){23}{\line(0,-1){.097826}}
\put(71,8.75){\makebox(0,0)[cc]{$p$}}
\put(124,39.25){\makebox(0,0)[cc]{$q$}}
\put(125.25,52.25){\makebox(0,0)[cc]{$\bbb_1$}}
\put(96.75,81.5){\makebox(0,0)[cc]{$\bbb_2$}}
\put(14.5,50.5){\makebox(0,0)[cc]{$\bbb_{N'}$}}
\put(45,25.5){\makebox(0,0)[cc]{$\Delta$}}
\put(71.25,35.5){\makebox(0,0)[cc]{$\Pi$}}
\end{picture}

\caption{$u$-close hub.} \label{fig3}
\end{figure}
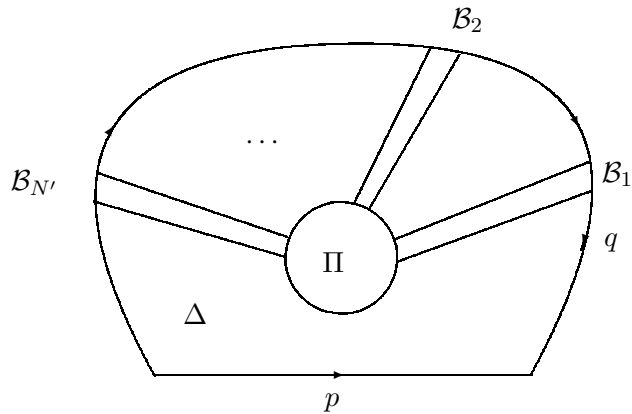

4. Suppose that there are two $k$-edges $e, e'$ in $q$ which are
connected by a $k$-band $\bbb$. Consider the subdiagram $\Delta'$
bounded by $\bbb$ and $q$ ($\Delta'$ does not contain $\bbb$). We
can assume that $\Delta'$ does not contain $k$-bands connecting
two edges on $q\cup\partial\Delta'$.  If there was a hub in
$\Delta'$ then by \cite[Lemmas 10.4 and 10.3]{OScol} there would
be a $u$-close hub there. Hence by part 1 there are no hubs in
$\Delta'$. Hence there are no $k$-edges between $e$ and $e'$ in
$q$, therefore $e$ and $e'$ belong to a subpath of $q$ labeled by
a cyclic shift $u_1$ of $u$. The word that labels the subpath of
$q$ starting at $e$ and ending at $e'$ is equal in $\hhh_1$ to the
word written on a side of $\bbb$ that is farther from $p$. Thus
the cyclic shift $u_1$ is equal in $\hhh_1$ to a word with fewer
$k$-letters, contrary the assumption that $u$ is minimal.

Thus there are no two $k$-edges in $q$ connected by a $k$-band.
Then, by part 3, the number of $k$-edges on $q$ is bounded by a
linear function in $|w|$. If $u$ contains a $k$-letter, we get a
recursive bound for $n$. Hence we may suppose that $k$-letters do
not occur in $u$.

5. The equality $w=u^n$ in $\hhh$ and the minimality of $w$
implies now that $w$ has no $k$-letters as well. Then, by
\cite[Lemmas 10.4, 10.3]{OScol}, $\Delta$ has no hubs, and by
\cite[Lemma 3.11]{OScol}, $\Delta$ has no $k$-annuli. So $\Delta$
has no $k$-edges at all. Therefore $w=u^n$ in $\hhh_1$. Now Lemma
3.11 \cite{OScol} implies that we may assume that $u$ contains no
$\theta$-letters. Indeed, again if $u$ contains a $\theta$-letter
and $n>|w|$, one of $\theta$-bands must connect two edges on the
subpath labeled by a cyclic shift of $u$, contrary the assumption
that $u$ is cyclically minimal. Since $w$ is minimal and $w=u^n$,
the word $w$ has no $\theta$-letters too, and there are no
$\theta$-edges in $\Delta$ by \cite[Lemma 3.11]{OScol}. Finally,
once again making use of the cyclic minimality of $u$ and the
minimality of $w$, we conclude that $u$ and $w$ have no
$a$-letters, and $\Delta$ has no $a$-bands for $a\notin\aaa(P_1)$.
Since every $x$-cell must be a member of an $a$-band for $a\notin
\aaa(P_1)$, all cells in $\Delta$ are $\gr$-cells. It remains to
use the solvability of the power problem in $\gr$ (and in its free
product with a free group generated by $x$-letters).
\endproof

{\em Proof of Theorem \ref{th1}.} Let $G$ be a group with solvable power problem. Using a
sequence of HNN extensions as in Lemma \ref{lm1}, we can embed $G$
into a group $G_1$ with solvable power problem
where every two elements of the same order are
conjugate. Thus the conjugacy problem in $G_1$ is decidable. By
Lemmas \ref{n3} and \ref{lm2}, $G_1$ can be embedded into a 2-generated group
$G_2$ with solvable power and conjugacy problems. Then
Lemma \ref{lm3} allows us to embed $G_2$ into a finitely presented group
$G_3$ with solvable power and conjugacy problems.
Finally applying Lemmas \ref{n3} and \ref{lm2}
again we embed $G_3$ into a $2$-generated finitely presented (by Lemma \ref{fp})
group $G_4$ with solvable power and conjugacy problems.
\medskip

Theorem \ref{th2} follows immediately from Lemmas \ref{n1}(6a), \ref{n3} and
\cite[Theorem 1.1]{OScol}.

\begin{minipage}[t]{2.9 in}
\noindent Alexander Yu. Ol'shanskii\\ Department of Mathematics\\
Vanderbilt University \\ alexander.olshanskiy@vanderbilt.edu\\
 and\\ Department of
Higher Algebra\\ MEHMAT\\
 Moscow State University\\
olshan@shabol.math.msu.su\\
\end{minipage}
\begin{minipage}[t]{2.6 in}
\noindent Mark V. Sapir\\ Department of Mathematics\\
Vanderbilt University\\
msapir@math.vanderbilt.edu\\
\end{minipage}

\end{document}